\documentclass[twoside]{article}
\usepackage{amsfonts}
\usepackage{tipa}
\usepackage{pifont}
\usepackage{amsmath}
\usepackage{amssymb}
\renewcommand\arraystretch{1.2}

\newcommand\heiti{\CJKfamily{hei}}
\newcommand\songti{\CJKfamily{song}}
\usepackage{CJK,indentfirst,amsmath,amsfonts,amssymb,amsthm,cite}
\begin{CJK*}{GBK}{song}
\usepackage{fancyhdr}
\usepackage[dvips]{graphics}
\usepackage[dvips]{color}

\begin{document}
\pagestyle{fancy}
\setcounter{page}{1} 
\newcounter{jie}

\renewenvironment{proof}[1][\heiti 证明]{\textbf{#1.} }{
\begin{flushright}
$\Box$
\end{flushright}
}
\fancyhf{} 
\fancyhead[RE]{\footnotesize{\songti  }} 
\fancyhead[LO]{\footnotesize{\songti  }} 
\fancyhead[LE,RO]{\thepage} 
\fancyhead[CE]{\footnotesize{\songti }} 
\fancyhead[CO]{\footnotesize{\songti  }} 
\renewcommand\headrulewidth{0.4pt} 

\thispagestyle{empty} \vspace*{-10mm} \noindent \hbox to \textwidth{
\renewcommand\arraystretch{1.0}
$\begin{array}{c} \mbox{{}}\\
\mbox{{\footnotesize {\bf}}}\\
\end{array}$
\hfill }
\vspace*{0.4cm}
\begin{center}
{
{\LARGE\bf On the number of solutions of the generalized Ramanujan-Nagell equation $D_{1}x^2+D_{2}^{m}=2^{n+2}$}}\\[0.6cm]
{\normalsize Li Jianghua\footnote{This work is supported by the
 N. S. F. of Shaanxi Province (2012K06-43) and Foundation of Shaanxi Educational Committee
 (12JK0874). E-mail: jianghuali@xaut.edu.cn.
}}\\[0.1cm]

\begin{minipage}{\textwidth}
\begin{minipage}[t]{15.5cm}
{}
\end{minipage}
\end{minipage}
College of Science, Xi'an University of Technology, Xi'an, Shaanxi,
P.R.China
\\[0.8cm]
\end{center}

{\small
\begin{minipage}{\textwidth}
\begin{minipage}[t]{13.2cm}
{\bf Abstract}  { Let $D_{1}$, $D_{2}$ be coprime odd integers with
min$(D_{1},D_{2})>1$, and let $N(D_{1},D_{2})$ denote the number of
positive integer solutions $(x, m, n)$ of the equation
$D_{1}x^2+D_{2}^{m}=2^{n+2}$. In this paper, we prove that
$N(D_{1},D_{2})\leq 2$ except for $N(3,5)=N(5,3)=4$ and
$N(13,3)=N(31,97)=3$.}
\end{minipage}

\ \

\begin{minipage}[t]{1.4cm}
 {\bf Keywords}
 \end{minipage} {\ \ \  Exponential diophantine equation; generalized Ramanujan-Nagell equation; nu\\-mber of solutions;
  upper bound.}


\end{minipage}
 }

\section*{\S 1. Introduction}
Let $\mathbb{Z}, \mathbb{N}$ be the sets of all integers and
positive integers respectively. Let $D_1$, $D_2$ be coprime positive
odd integers with $D_2>1$. In 1913, S. Ramanujan [18] conjectured
that all the solutions $(x,n)$ of the equation
$$
x^2+7=2^{n+2},\ x, n \in\mathbb{N}
$$
are given by $(x,n)=(1,1),\ (3,2),\ (5,3),\ (11,5)$ and $(181,13)$.
Afterwards, W. Ljunggren [11] posed the same problem and T. Nagell
[17] solved it in 1948. Subsequently, the equation
$$D_{1}x^2+D_{2}=2^{n+2},\ x, n \in\mathbb{N} \eqno{(1.1)}$$ is usually called the generalized Ramanujan-Nagell equation,
which was solved by Y. Bugeaud and T. N. Shorey [7].

In this paper we deal with the number of solutions $(x,m,n)$ of the
equation  $$D_{1}x^2+D_{2}^{m}=2^{n+2},\ x, m, n \in\mathbb{N},
\eqno{(1.2)}$$ which is an exponential extension of (1.1). Let
$N(D_{1},D_{2})$ denote the number of solutions $(x,m,n)$ of (1.2).
For $D_{1}=1$, sum up the results of [4] and [9], we have:

{\bf Theorem A.} If $D_1=1$, then $N(1,D_2)\leq 1$ except for
$N(1,7)=6$, $N(1,23)=2$ and $N(1,2^r-1)=2$, where $r\in \mathbb{N}$
with $r>3$.

For $D_1>1$, we prove a general result as follows:

{\bf Theorem B.} If $D_1>1$, then $N(D_1,D_2)\leq 2$ except for
$N(3,5)=N(5,3)=4$ and $N(13,3)=N(31,97)=3$.

\section*{\S 2. Preliminaries}

{\bf Lemma 2.1.}([10, Formula 1.76]) For any positive integer $k$
and any complex numbers $\alpha,\ \beta$, we have
$$\alpha^k+\beta^k=\sum_{i=0}^{[k/2]}(-1)^i\left[k \atop i \right](\alpha+\beta)^{k-2i}(\alpha\beta)^i,$$
where $[k/2]$ is the integral part of $k/2$,
$$\left[k \atop i \right]=\frac{(k-i-1)!k}{(k-2i)!i!} \in \mathbb{N},\ i=0, 1,
..., [k/2].$$

{\bf Lemma 2.2.}([3]) Let $p$ be an odd prime, and let $X$ be an
integer with $|X|>1$. If $q$ is a prime divisor of
$(x^p-1)/(x-1)$, then either $q=p$ or $q\equiv 1 ({\rm mod}\ 2p)$.
Further, if $p\mid(x^p-1)/(x-1)$, then $p\parallel(x^p-1)/(x-1)$.

For any nonnegative integer $k$, let $F_k$ and $L_k$ denote the
$k$-th Fibonacci number and Lucas number respectively.

{\bf Lemma 2.3.}([14, pp. 60-61])

(i) $2\mid F_kL_k$ if and only if $3\mid k$.

(ii) \[\rm {gcd}(F_k, L_k)=\begin{cases}1, &{\rm if} \ 3\nmid k,\\
2, &{\rm if}\ 3\mid k.\end{cases}\]

(iii) $L_k^2-5F_k^2=(-1)^k4$.

(iv) Every solution $(u,v)$ of the equation $$u^2-5v^2 =\pm 4,\
 u,\ v \in \mathbb{N}$$
 can be expressed as $(u,v)=(L_k,F_k)$, where $k\in \mathbb{N}.$

{\bf Lemma 2.4.} ([6]) The equation
$$F_k =z^n,\ k,\ z, \ n \in \mathbb{N},\ z>1,\ n>1$$ has only the
solutions $(k,z,n)=(6,2,3)$ and $(12,12,2)$. The equation
$$L_k=z^n,\ k,\ z,\ n\in \mathbb{N},\ z>1,\ n>1$$ has only the
solution $(k,z,n)=(3,2,2).$

{\bf Lemma 2.5.} ([16]) The equation
$$x^3+1 =3y^2,\ x,\ y \in \mathbb{N}$$ has no solution $(x,y)$.

{\bf Lemma 2.6.} ([15]) Let $p$ be an odd prime. The equation
$$x^2+x+1=3y^p,\ x,\ y,  \in \mathbb{Z},\ \mid x\mid>1,\ y>1$$ has no
solution $(x,y)$.

{\bf Lemma 2.7.} ([5]) The equation
$$\frac{2^r+1}{2+1}=y^n,\ r,\ y, \ n \in \mathbb{N},\ y>1,\ n>1$$ has no
solution $(r,y,n)$.

{\bf Lemma 2.8.} ([13]) The equation
$$x^m-y^n=1, \ x,\ y, \ m, \ n \in \mathbb{N},\ {\rm min}(x,y,m,n)>1$$ has
only the solution $(x,y,m,n)=(3,2,2,3)$.

{\bf Lemma 2.9.} The equation
$$2^r+1=3^sy^n, \ r,\ s, \ y, \ n \in \mathbb{N},\ 3\nmid y,\ y>1,\ n>1\eqno{(2.1)}$$ has
no solution $(r,s,y,n)$.

{\bf Proof.} By Lemma 2.7, (2.1) has no solution $(r,s,y,n)$ with
$s=1$.

If $s>1$ and $2\mid n$, then we have $2\nmid r$, $3\mid r$ and
$(2^{r/3})^3+1=3^s(y^{n/2})^2$ by (2.1). But, since $y>1$, by Lemmas
2.5 and 2.8, it is impossible.

If $s>1$ and $2\nmid n$, then $2\nmid r$, $3\mid r$ and $n$ has an
odd prime divisor $p$. Since $2^r+1=(2^{r/3}+1)(2^{2r/3}-2^{r/3}+1)$
and gcd$(2^{r/3}+1,\ 2^{2r/3}-2^{r/3}+1)=3$, we get from (2.1) that
$2^{r/3}+1=3^{s-1}a^p$ and$$ 2^{2r/3}-2^{r/3}+1=3b^p \eqno{(2.2)}$$
where $a,b\in \mathbb{N}$ with $ab=y^{n/p}$. But, since $2^{r/3}>1$,
by Lemma 2.6, (2.2) is impossible. Thus, the lemma is proved.

{\bf Lemma 2.10.}  The equation
$$2^r-1=3^sy^n, \ r,\ s, \ y, \ n \in \mathbb{N},\ 3\nmid y,\ y>1,\ n>1\eqno{(2.3)}$$ has
no solution $(r,s,y,n)$.

{\bf Proof.} We see from (2.3) that $r$ must be even. Since
gcd$(2^{r/2}+1, 2^{r/2}-1)=1$, by (2.3), we have

\[2^{r/2}+1=\begin{cases}3^sa^n, \\ b^n,\end{cases}2^{r/2}-1=\begin{cases}b^n, \\ 3^sa^n,\end{cases}\ y=ab,\ a,b\in \mathbb{N},
\eqno{(2.4)}\] However, since $y>1$ and $3\nmid y$, by Lemma 2.8,
(2.4) is impossible. Thus, the lemma is proved.

{\bf Lemma 2.11.}  The equation
$$2^r\cdot3^s+1=y^n, \ r,\ s, \ y, \ n \in \mathbb{N},\ y>1,\ n>1\eqno{(2.5)}$$ has
only the solutions $(r,s,y,n)=(3,1,5,2),\ (4,1,7,2)$ and
$(5,2,17,2)$.

{\bf Proof.} If $2\mid n$, since $2\nmid y$, then we have $r\geq3$
and gcd$(y^{n/2}+1, y^{n/2}-1)=2$. Hence, by (2.5), we get

\[y^{n/2}+1=\begin{cases}2^{r-1}, \\ 2\cdot3^s,\end{cases}y^{n/2}-1=\begin{cases}2\cdot3^s, \\ 2^{r-1},\end{cases}\eqno{(2.6)}\]
whence we obtain  \[1=\begin{cases}2^{r-2}-3^s, \\
3^s-2^{r-2}.\end{cases}\eqno{(2.7)}.\] Apply Lemma 2.8 to (2.7), we
get $(r,s,y,n)=(3,1,5,2),\ (4,1,7,2)$ and $(5,2,17,2)$ by (2.6).

If $2\nmid n$, since $n>1$, then $n$ has an odd prime divisor $p$.
By (2.5), we get$$2^r\cdot3^s=y^n-1=(z-1)(z^{p-1}+z^{p-2}+...+1),\
 z=y^{n/p}. \eqno{(2.8)}$$
Since $2\nmid z^{p-1}+z^{p-2}+...+1$ and $3\not\equiv1({\rm mod} \
2p)$, by Lemma 2.2, we see from (2.8) that $3\geq
z^{p-1}+z^{p-2}+...+1$, a contradiction. Thus, the lemma is
proved.

{\bf Lemma 2.12.}  The equation
$$2^r\cdot3^s-1=y^n, \ r,\ s, \ y, \ n \in \mathbb{N},\ y>1,\ n>1\eqno{(2.9)}$$ has
no solution $(r,s,y,n)$.

{\bf Proof.} Since $(-1/3)=-1$, where $(\ast/\ast)$ is the Jacobi
symbol, we see from (2.9) that $n$ must be odd. Since $n>1$, $n$ has
an odd prime divisor $p$, and by (2.9), we have
$$2^r\cdot3^s=y^n+1=(z+1)(z^{p-1}-z^{p-2}+...+1),\
 z=y^{n/p}. \eqno{(2.10)}$$
Apply Lemma 2.2 to (2.10), we get $3\geq z^{p-1}-z^{p-2}+...+1$, a
contradiction. Thus, the lemma is proved.

{\bf Lemma 2.13.}  ([7, Lemma 1]) If the equation
$$D_{1}X^2+D_{2}Y^{2}=2^{Z+2}, \ X,\ Y, \ Z \in \mathbb{Z},\ {\rm gcd}(X,Y)=1,\ Z>0\eqno{(2.11)}$$
has solutions $(X,Y,Z)$, then it has a unique positive integer
solution $(X_1,Y_1,Z_1)$ satisfying $Z_1\leq Z$, where $Z$ through
all solutions $(X,Y,Z)$ of (2.11). Such $(X_1,Y_1,Z_1)$ is called
the least solution of (2.11). Every solution  $(X,Y,Z)$ of (2.11)
can be expressed as
$$Z=Z_1t,\ \ t\in  \mathbb{N},\ 2\nmid t \ {\rm if}\ D_{1}>1,$$
$$\frac{X\sqrt{D_1}+Y\sqrt{-D_2}}{2}=\lambda_1\left(\frac{X_1\sqrt{D_1}+\lambda_2Y_1\sqrt{-D_2}}{2}\right)^t,
\ \lambda_1,\ \lambda_2 \in \{\pm 1\}.$$

By Lemma 2.13, we can obtain the following lemma  immediately.

{\bf Lemma 2.14.} If $(X,Y,Z)$ and $(X',Y',Z')$ are two positive
integer solutions of (2.11) with $(X,Y,Z)\neq (X',Y',Z')$, then
$Z\neq Z'$.

{\bf Lemma 2.15.}([8, Lemma 3]) Let min$(D_1,D_2)>1$ and $D=D_1D_2$.
If (2.11) has solutions $(X,Y,Z)$, then the equation
$$X'^2+DY'^{2}=2^{Z'+2}, \ X',\ Y', \ Z' \in \mathbb{Z},\ {\rm gcd}(X',Y')=1,\ Z'>0\eqno{(2.12)}$$
has solutions $(X',Y',Z')$. Moreover, its least solution
$(X_1',Y_1',Z_1')$ satisfies $X'_1=\frac{1}{2}\mid
D_{1}X_1^2-D_{2}Y_1^{2}\mid$, $Y'_1=X_1Y_1$ and $Z'_1=2Z_1$, where
$(X_1,Y_1,Z_1)$ is the least solution of (2.11).

{\bf Lemma 2.16.} Let min$(D_1,D_2)>1$ and $D=D_1D_2$. If (2.1) has
a solution $(X,Y,Z)$, then (2.2) has no solution $(X',Y',Z')$ with
$Z'=Z$.

{\bf Proof.} By Lemma 2.13, if $(X,Y,Z)$ and $(X',Y',Z')$ are
solutions of (2.11) and (2.12), then we have
$$Z=Z_1t,\ t\in \mathbb{N},\ 2\nmid t \eqno{(2.13)}$$
and$$Z'=Z'_1t',\ t'\in \mathbb{N},  \eqno{(2.14)}$$ where
$(X_1,Y_1,Z_1)$ and $(X_1',Y_1',Z_1')$ are least solutions of (2.11)
and (2.12) respectively. Further, by Lemma 2.15, we have
$Z'_1=2Z_1$. Substituting it into (2.14), we get $Z'=2Z_1t'$. Since
$2\nmid t$, we obtain $Z'\neq Z$ by (2.13). Thus, the lemma is
proved.

{\bf Lemma 2.17.} Let $D$ be a positive integer. Further, let (2.12)
have solutions  $(X',Y',Z')$ and $(X_1',Y_1',Z_1')$ is its least
solution. If $(y,z)$ is a solution of the equation$$1+Dy^2=2^{z+2},\
y,\ z\in\mathbb{N}, \eqno(2.15)$$ then $X'_1=1$ and
$(y,z)=(Y_1',Z_1')$ except for $D=7$ and $(y,z)=(3,4).$

{\bf Proof.} Under the assumption, (2.12) has the solution
$(X',Y',Z')=(1,y,z)$. By Lemma 2.13, we have
$$z=Z'_1t,\ t\in\mathbb{N},\eqno{(2.16)}$$
$$\frac{1+y\sqrt{-D}}{2}=\lambda_1\left(\frac{X'_1+\lambda_2Y'_1\sqrt{-D}}{2}\right)^t,
\ \lambda_1,\ \lambda_2 \in \{\pm 1\}.\eqno{(2.17)}$$

If $2\mid t$,
let$$\frac{a+b\sqrt{-D}}{2}=\left(\frac{X'_1+\lambda_2Y'_1\sqrt{-D}}{2}\right)^{t/2}.\eqno{(2.18)}$$
By Lemma 2.13, then $a,\ b$ are integers satisfying
$$a^2+Db^2=2^{Z'_1t/2+2}=2^{z/2+2},\ {\rm gcd}(a,b)=1.\eqno{(2.19)}$$
Substituting (2.18) into (2.17), we get
$$a^2-Db^2=2\lambda_1,\ y=ab\lambda_1,
\ \lambda_1\in \{\pm 1\}.\eqno{(2.20)}$$ The combination of (2.19)
and the first equality of (2.20) yields $\lambda_1=1$ and
$a^2=2^{z/2+1}+1$. Hence, by Lemma 2.8, we get $a=3,\ z=4,\ D=7,\
b=1$ and $y=3$ by (2.20).

If $2\nmid t$ and $t>1$, let
$$\alpha=\frac{X'_1+Y'_1\sqrt{-D}}{2},\
\beta=\frac{X'_1-Y'_1\sqrt{-D}}{2},\eqno{(2.21)}$$then from (2.17)
we get $$1=\lambda_1(\alpha^t+\beta^t).\eqno{(2.22)}$$ Since
$\alpha+\beta=X'_1$ and $\alpha\beta=2^{Z'_1}$ by (2.21), apply
Lemma 2.1 to (2.22) we have
$$1=\lambda_1\sum_{i=0}^{(t-1)/2}(-1)^i\left[t \atop i \right](\alpha+\beta)^{t-2i}(\alpha\beta)^i$$
$$=\lambda_1X'_1\sum_{i=0}^{(t-1)/2}(-1)^i\left[t \atop i \right]X'^{t-2i-1}_1 2^{Z'_1i}.\eqno{(2.23)}$$
From (2.23), we obtain $X'_1=1$, $\lambda_1=1$ and
$$t=\left[t \atop 1 \right]=2^{Z'_1}\sum_{j=2}^{(t-1)/2}(-1)^j\left[t \atop j \right]2^{Z'_1(j-2)},$$
a contradiction. It implies that $t=1$ if $2\nmid t$. Thus, by
(2.16) and (2.17), we get $X'_1=1$ and $(y,z)=(Y'_1,Z'_1)$. The
lemma is proved.

Let $\alpha,\ \beta$ be algebraic integers. If $(\alpha+\beta)^2$
and $\alpha\beta$ are nonzero coprime integer and $\alpha/\beta$ is
not a root of unity, then $(\alpha,\ \beta)$ is called a Lehmer
pair. Let $a=(\alpha+\beta)^2$ and $b=\alpha\beta$. Then we have
$$\alpha=\frac{1}{2}(\sqrt{a}+\lambda\sqrt{c}),\ \beta=\frac{1}{2}(\sqrt{a}-\lambda\sqrt{c}),\ \lambda\in \{-1,1\},$$
where $c=a^2-4b$. The pair $(a,c)$ is called the parameter of the
Lehmer pair $(\alpha,\beta)$. Two Lehmer pairs $(\alpha_1,\beta_1)$
and $(\alpha_2,\beta_2)$ are equivalent if
$\alpha_1/\alpha_2=\beta_1/\beta_2\in \{\pm1,\pm\sqrt{-1}\}.$ Given
a Lehmer pair $(\alpha,\beta)$, one defines the corresponding Lehmer
numbers by
\[L_k(\alpha,\beta)=\begin{cases}\displaystyle \frac{\alpha^k-\beta^k}{\alpha-\beta}, &2\nmid k,\\
\displaystyle\frac{\alpha^k-\beta^k}{\alpha^2-\beta^2}, &2\mid k,
\end{cases}\ \  k\in\mathbb{N}.\eqno{(2.24)}\]
Lehmer numbers are nonzero integers. For equivalent Lehmer pairs
$(\alpha_1,\beta_1)$ and $(\alpha_2,\beta_2)$, we have
$L_k(\alpha_1,\beta_1)=\pm L_k(\alpha_2,\beta_2)$ $(k\in
\mathbb{N})$. A prime $q$ is called a primitive divisor of
$L_k(\alpha,\beta)$ $(k>1)$ if $p\mid L_k(\alpha,\beta)$ and $p\nmid
(\alpha^2-\beta^2)^2L_1(\alpha,\beta)...L_{k-1}(\alpha,\beta)$. A
Lehmer pair $(\alpha,\beta)$ such that $L_k(\alpha,\beta)$ has no
primitive divisor will be called a $k$-defective Lehmer pair.
Further, a positive integer $k$ is called totally non-defective if
no Lehmer pair is $k$-defective.

{\bf Lemma 2.18.} ([1],[9]) Let $k$ satisfy $6<k\leq 30$ and
$2\nmid k$. Then, up to equivalence, all parameters of
$k$-defective Lehmer pairs are given as follows:

(i) $k=7$, $(a,c)=(1,-7),\ (1,-19),\ (3,-5),\ (5,-7),\ (13,-3),\
(14,-22).$

(ii) $k=9$, $(a,c)=(5,-3),\ (7,-1),\ (7,-5).$

(iii) $k=13$, $(a,c)=(1,-7).$

(iv) $k=15$, $(a,c)=(7,-1),\ (10,-2).$

{\bf Lemma 2.19.} ([2], Theorem 1.4) If $k>30$, then $k$ is totally
non-defective.

\section*{\S 3. Further lemmas on the solutions of (1.2)}

Let $D_1>1$. We first consider the solutions $(x,m,n)$ of (1.2) with
$2\nmid m$. Then (2.11) has the solution $$(X,Y,Z)=(x,\
D_2^{(m-1)/2},n).\eqno{(3.1)}$$Since min$(D_1,D_2)>1$, apply Lemma
2.13 to (3.1), we get $$n=Z_1t,\ t\in\mathbb{N},\ 2\nmid
t,\eqno(3.2)$$
$$\frac{x\sqrt{D_1}+D_2^{(m-1)/2}\sqrt{-D_2}}{2}=\lambda_1\left(\frac{X_1\sqrt{D_1}+\lambda_2Y_1\sqrt{-D_2}}{2}\right)^t,
\ \lambda_1,\ \lambda_2 \in \{\pm 1\}.\eqno{(3.3)}$$ where
$(X_1,Y_1,Z_1)$ is the least solution of (2.11). Let
$$\alpha=\frac{X_1\sqrt{D_1}+Y_1\sqrt{-D_2}}{2},\ \beta=\frac{X_1\sqrt{D_1}-Y_1\sqrt{-D_2}}{2}.\eqno{(3.4)}$$
Since $X_1,\ Y_1$ and $Z_1$ are positive integers satisfying
$$D_{1}X_1^2+D_{2}Y_1^{2}=2^{Z_1+2}, \ {\rm gcd}(D_{1}X_1^2,D_{2}Y_1^{2})=1,\eqno{(3.5)}$$
$\alpha$ and $\beta$ are roots of
$z^4-\frac{1}{2}(D_{1}X_1^2-D_{2}Y_1^{2})z^2+2^{2Z_1}=0.$ Notice
that $(\alpha+\beta)^2=D_{1}X_1^2$ and $\alpha\beta=2^{Z_1}$ are
coprime positive integers, and
$\alpha/\beta=(\frac{1}{2}(D_1X_1^2-D_2Y_1^2)+X_1Y_1\sqrt{-D_1D_2})/2^{Z_1}$is
not a root of unity. Then $(\alpha,\beta)$ is a Lehmer pair with
parameter $(D_1X_1^2, -D_2Y_1^2)$. Let $L_k(\alpha,\beta)$ $(k\in
\mathbb{N})$ denote the corresponding Lehmer numbers defined as in
(2.24). By (3.3) and (3.4), we have
$$D_2^{(m-1)/2}=Y_1\mid L_t(\alpha,\beta)\mid.\eqno{(3.6)}$$
Since $(\alpha^2-\beta^2)^2=-D_1D_2X_1^2Y_1^2$, we find from (3.6)
that the Lehmer number $L_t(\alpha,\beta)$ has no primitive divisor.
Therefore, we have the following result.

{\bf Lemma 3.1.} If $D_1>1$, then all solutions $(x,m,n)$ of (1.2)
with $2\nmid m$ are given as follows:

(i) $t=9$, $(D_1,D_2)=(5,3),\ (x,m,n)=(19,5,9).$

(ii) $t=7$, $(D_1,D_2)=(3,5),\ (x,m,n)=(13,1,7).$

(iii) $t=7$, $(D_1,D_2)=(13,3),\ (x,m,n)=(71,1,14).$

(iv) $t=5$, $(D_1,D_2)=(5,3),\ (x,m,n)=(5,1,5).$

(v) $t=5$, $(D_1,D_2)=(21,11),\ (x,m,n)=(79,1,15).$

(vi) $t=5$, $(D_1,D_2)=(3,29),\ (x,m,n)=(209,1,15).$

(vii) $t=5$, $(D_1,D_2)=(3,5),\ (x,m,n)=(1,3,5).$

(viii) $t=5$, $(D_1,D_2)=(11,5),\ (x,m,n)=(19,3,10).$

(ix) $t=3$, $$D_1X_1^2=2^{Z_1}-\lambda,\ D_2=3\cdot2^{Z_1}+\lambda,\
\lambda\in
\{\pm1\},\eqno{(3.7)}$$$(x,m,n)=(X_1(2^{Z_1+1}-\lambda),1,3Z_1).$

(x) $t=3$, $(D_1,D_2)=(5,3),\ (x,m,n)=(1,3,3).$

(xi) $t=3$, $(D_1,D_2)=(13,3),\ (x,m,n)=(1,5,6).$

(xii) $t=1$, $Y_1=D_2^{(m-1)/2},\ (x,m,n)=(X_1,m,Z_1).$

{\bf Proof.} By Lemma 2.19, we have $t\leq30.$ Further, since
$2\nmid t$, by Lemma 2.18, if $7\leq t\leq 30$, then (1.2) has only
the solutions (i), (ii) and (iii) satisfying $2\nmid m.$

For $t=5$, apply Lemma 2.1 to (3.6), we get
$$D_2^{(m-1)/2}=Y_1\mid (D_2Y_1^2)^2-5\cdot 2^{Z_1}(D_1Y_1^2)+5\cdot 2^{2Z_1}\mid.\eqno{(3.8)}$$

If $m=1$, then from (3.8) we get $Y_1=1$ and

$$D_2^2-5\cdot2^{Z_1}D_2+5\cdot2^{2Z_1}=(D_2-5\cdot2^{Z_1-1})^2-5\cdot2^{2Z_1-2}=\lambda,\ \lambda\in\{\pm1\}.\eqno{(3.9)}$$

When $Z_1=1$, by (3.9), we have $\mid D_2-5\cdot2^{Z_1-1}\mid=\mid
D_2-5\mid=2.$ It implies that $D_2=3$ or 7, and by (3.5), $X_1=1$
and $D_1=5$ or 1. Since $D_1>1,$ we get $(D_1,D_2)=(5,3)$ and the
solution (iv).

When $Z_1>1$, by (3.9), we get
$$(D_2-5\cdot2^{Z_1-1})^2-5(2^{Z_1-1})^2=1.\eqno{(3.10)}$$
Apply Lemma 2.3 to (3.10), we have $$\mid
D_2-5\cdot2^{Z_1-1}\mid=\frac{1}{2}L_{6I+6},\
2^{Z_1-1}=\frac{1}{2}F_{6I+6}, \ I\in \mathbb{Z},\ I\geq
0.\eqno{(3.11)}$$ Further, by Lemma 2.4, we see from the second
equality of (3.11) that $I=0$ and $Z_1=3$. Hence, by the first
equality of (3.11), we get $D_2=11$ or 29. Further, by (3.5), we
have $X_1=1$ and $D_1=21$ or 3. Thus, by (3.3), the solutions (v)
and (vi) are obtained.

If $m>1$ and $5\nmid D_2$, since gcd$(D_2, 5\cdot2^{Z_1})=1,$ then
from (3.8) we get $Y_1=D_2^{(m-1)/2}$ and
$$(D_2^m-5\cdot2^{Z_1-1})^2-5(2^{Z_1-1})^2=\lambda,\
\lambda\in\{\pm1\}.\eqno{(3.12)}$$ Since $m>1$ and min$(D_1,D_2)>1$,
using the same method as in the case $m=1$, we can prove that (3.12)
is impossible.

If $m>1$ and $5\mid D_2$, then we have

$$Y_1=\frac{1}{5}{D_2^{(m-1)/2}} \eqno{(3.13)}$$
and
$$(2^{Z_1-1})^2-5\left(\frac{1}{125}D_2^m-2^{Z_1-1} \right)^2=\lambda, \ \lambda \in \{\pm1\}. \eqno{(3.14)}$$

When $Z_1=1$, by (3.14), we get $\frac{D_2^m}{125}-1=0$, and hence,
we have $D_2=5$ and $m=3$. It implies that $(D_1,D_2)=(3,5)$ and the
solution (vii) is obtained.

 When $Z_1=2$, by (3.14), we have $|\frac{1}{125}{D_2^m}-2|=1$, whence we get $D_2=5$ and $m=3$. Further, by (3.3), (3.5) and (3.13), we obtain $X_1=Y_1=1$, $D_1=11$ and the solution (viii).

  For $t=3$, by (3.8), we have $$D_2^{(m-1)/2}=Y_1|D_2Y_1^2-3\cdot2^{Z_1}|.  \eqno {(3.15)}$$

  If $m=1$, then from (3.15) we get $Y_1=1$ and $D_2=3\cdot2^{Z_1}+\lambda$, where $\lambda \in \{\pm1\}$. Hence, by (3.5), $D_1$, $D_2$ satisfy (3.7) and the solution (ix) is obtained.

If $m>1$ and $3\nmid D_2$, then we have $Y_1=D_2^{(m-1)/2}$ and
$${D_2^m}=3\cdot2^{Z_1}+\lambda, \ \lambda \in \{\pm1\}. \eqno
{(3.16)}$$ However, since $D_2>1, m>1$ and $2\nmid m$, by Lemmas
2.11 and 2.12, (3.18) is impossible.

If $m>1$ and $3\mid D_2$, then we have
$$Y_1=\frac{1}{3}D_2^{(m-1)/2} \eqno{(3.17)}$$ and
$$\frac{1}{27}{D_2^m}-2^{Z_1}=\lambda,\  \lambda \in \{\pm1\}.
\eqno{(3.18)}$$ Let $D_2=3^ID$, where $I,D \in \mathbb{N}$ with
$3\nmid D$. Then (3.18) can be written as
$$3^{Im-3}D^m-2^{Z_1}=\lambda,\  \lambda \in \{\pm1\}.
\eqno{(3.19)}$$

When $D=1$, by (3.19), we have $$3^{Im-3}-2^{Z_1}=\lambda,\  \lambda
\in \{\pm1\}. \eqno{(3.20)}$$ Since $m>1$ and $2\nmid m$, apply
Lemma 2.8 to (3.20), we get $(I,m,Z_1,\lambda)=(1,3,1,-1)$ and
$(1,5,3,1)$. Therefore, by (3.3), (3.5) and (3.17), the solutions
(x) and (xi) are obtained.

When $D>1$, by Lemmas 2.9 and 2.10, (3.19) is impossible.

For $t=1$, by (3.2) and (3.3), the solutions (xii) is obtained. To
sum up, the lemma is proved.

Let $N_1(D_1,D_2)$ denote the number of solutions $(x,m,n)$ of (1.2)
with $2\nmid m$. By Lemma 3.1, we can obtain the following lemma
immediately.

{\bf Lemma 3.2.} If $D_1>1$, then $N_1(D_1,D_2)\leq 1$ except for
the following cases:

(i) $N_1(3,5)=4, (x,m,n)=(1,1,1),(3,1,3),(1,3,5)$ and $(13,1,7)$.

(ii)  $N_1(5,3)=3, (x,m,n)=(1,1,1),(1,3,3),(5,1,5)$ and $(19,5,9)$.

(iii)  $N_1(13,3)=3, (x,m,n)=(1,1,2),(1,5,6)$ and $(71,1,14)$.

(iv) $N_1(11,5)=2, (x,m,n)=(1,1,2)$ and $(19,3,10)$.

(v) $N_1(21,11)=2, (x,m,n)=(1,1,3)$ and $(79,1,15)$.

(vi) $N_1(3,29)=2, (x,m,n)=(1,1,3)$ and $(209,1,15)$.

(v) If $D_1$ and $D_2$ satisfy (3.7) with $Z_1\geq 2$, then
$N_1(D_1,D_2)=2$, $(x,m,n)=(X_1,1,Z_1)$ and
$(X_1(2^{Z_1+1}-\lambda),1,3Z_1)$.

We next consider the solutions $(x,m,n)$ of (1.2) with $2|m$. Then
the equation
$$D_1{X^\prime}^2+D_2^2{Y^\prime}^2=2^{Z^\prime+2},X^\prime,Y^\prime,Z^\prime
\in \mathbb{Z}, \textrm{gcd}(X^\prime,Y^\prime)=1,Z^\prime>0
\eqno{(3.21)}$$ has the solution
$$(X^\prime,Y^\prime,Z^\prime)=(x,D_2^{(m-2)/2},n).\eqno{(3.22)}$$
Since $\textrm{min}(D_1,D_2)>1$, apply Lemma 2.13 to (3.22), we have
$$n=Z_1^\prime t^\prime,t^\prime \in \mathbb{N}, 2\nmid t^\prime,
\eqno{(3.23)}$$
$$\frac{x\sqrt{D_1}+D_2^{(m-2)/2}\sqrt{-D_2^2}}{2}=\lambda_1\left( \frac{X_1^\prime \sqrt{D_1}+\lambda_2Y_1^\prime\sqrt{-D_2^2} }{2}    \right)^{t^\prime},\ \lambda_1,\lambda_2 \in \{\pm 1\}, \eqno{(3.24)}$$
where $(X_1^\prime,Y_1^\prime,Z_1^\prime)$ is the least solution of
(3.21). Let
$$\alpha^\prime=\frac{X_1^\prime\sqrt{D_1}+Y_1^\prime\sqrt{-D_2^2}}{2},
\ \beta^\prime=\frac{X_1\sqrt{D_1}-Y_1^\prime\sqrt{-D_2^2}}{2}.
\eqno{(3.25)}$$ Then $(\alpha^\prime, \beta^\prime)$ is a Lehmer
pair with parameter $(D_1{X_1^\prime}^2,-D_2^2{Y_1^\prime}^2)$.
Further, let $L_k(\alpha^\prime, \beta^\prime) (k \in \mathbb{N})$
denote the corresponding Lehmer numbers. Form (3.24) and (3.25), we
have
$$D_2^{(m-2)/2}=Y_1^\prime|L_{t^\prime}(\alpha^\prime,
\beta^\prime)|. \eqno{(3.26)}$$ Since
$({\alpha^\prime}^2-{\beta^\prime}^2)^2=-D_1D_2^2{X_1^\prime}^2{Y_1^\prime}^2$,
we see from (3.26) that the Lehmer number
$L_{t^\prime}(\alpha^\prime, \beta^\prime)$ has no primitive
divisor. Therefore, using the same method as in the proof of Lemma
3.1, we can obtain the following lemma.

{\bf Lemma 3.3.} If $D_1>1$, then all the solutions $(x,m,n)$ of
(1.2) with $2|m$ are given as follows:

(i) $t^\prime=3, (D_1,D_2)=(7,3),(x,m,n)=(5,4,6)$.

(ii) $t^\prime=3, (D_1,D_2)=(7,5),(x,m,n)=(17,2,9)$.

(iii) $t^\prime=3, (D_1,D_2)=(15,7),(x,m,n)=(33,2,12)$.

(iv) $t^\prime=1,
Y_1^\prime={D_2}^{(m-2)/2},(x,m,n)=(X_1^\prime,m,Z_1^\prime)$.

Let $N_2(D_1,D_2)$ denote the number of solutions $(x,m,n)$ of (1.2)
with $2|m$. By Lemma 3.3, we have:

{\bf Lemma 3.4.} If $D_1>1$, then $N_2(D_1,D_2)\leq 1$ except for
the following cases:

(i) $N_2(7,3)=2, (x,m,n)=(1,2,2)$ and $(5,4,6)$.

(ii) $N_2(7,5)=2, (x,m,n)=(1,2,3)$ and $(17,2,9)$.

(iii) $N_2(15,7)=2, (x,m,n)=(1,2,4)$ and $(33,2,12)$.

\section*{\S 4. Proof of Theorem B}

{\bf Lemma 4.1.}([12]) The equation $$\frac{x^n+1}{x+1}=y^2,x,y,z
\in \mathbb{N}, x>1,n>1$$ has no solution $(x,y,n)$.

{\bf Lemma 4.2.} The equation $$2^{2r-3}-2^r+1=97^s,r,s\in
\mathbb{N}, r\geq 5 \eqno{(4.1)}$$ has only the solution
$(r,s)=(5,1)$.

{\bf Proof.} By (4.1), we have
$$2(2^{r-2}-1)^2=97^s+1.\eqno{(4.2)}$$ Since $r\geq 5$ and
$2^{r-2}-1$ has an odd prime divisor $p$ with $p\equiv 3 (\rm mod\
4)$, if $2|s$, then (4.2) is impossible. So we have $2\nmid s$.

Since $2\nmid s$ and $97+1=2.7^2$, we see from (4.2) that
$7|2^{r-2}-1$ and $$\frac{97^s+1}{97+1}={\left(
\frac{2^{r-2}-1}{7}\right)}^2. \eqno{(4.3)}$$ Apply Lemma 4.1 to
(4.3), we get $s=1$ and $r=5$. Thus, the lemma is proved.

{\bf Lemma 4.3.} The equation $$7x^2+25^{2y}=2^{z+2},x,y,z\in
\mathbb{N} \eqno{(4.4)}$$ has no solution $(x,y,z)$.

{\bf Proof.} We suppose that (4.4) has a solution $(x,y,z)$. Then
the equation  $$7X^2+25Y^2=2^{Z+2},X,Y,Z\in \mathbb{Z},
\textrm{gcd}(X,Y)=1,Z>0  \eqno{(4.5)}$$ has the solution
$$(X,Y,Z)=(x,5^{2y-1},z). \eqno{(4.6)}$$ Since
$(X_1,Y_1,Z_1)=(1,1,3)$ is the least solution of (4.5), apply
Lemma 2.13 to (4.6), we have $$z=3t,t\in \mathbb{N}, 2\nmid t,
t>1, \eqno{(4.7)}$$
$$\frac{x\sqrt{7}+5^{2k-1}\sqrt{-25}}{2}=\lambda_1{\left(\frac{\sqrt{7}+\lambda_2\sqrt{-25}}{2}      \right)}^{t},\lambda_1,\lambda_2 \in \{\pm1\}.\eqno{(4.8)}$$

Let $$\alpha=\frac{\sqrt{7}+\sqrt{-25}}{2},\ \
\beta=\frac{\sqrt{7}-\sqrt{-25}}{2}. \eqno{(4.9)}$$ Then
$(\alpha,\beta)$ is a Lehmer pair with parameter $(7,-25)$.
Further, let $L_k(\alpha,\beta) (k\in \mathbb{N})$ denote the
corresponding Lehmer numbers. By (4.8) and (4.9), we have
$$5^{2k-1}=|L_t(\alpha,\beta)|. \eqno{(4.10)}$$ Since
$(\alpha^2-\beta^2)=-7\cdot25$, we see from (4.10) that the Lehmer
number $L_t(\alpha,\beta)$ has no primitive divisor. Therefore, by
Lemmas 2.18 and 2.19, we get $t\leq 5$. However, since
$L_5(\alpha,\beta)=-55$ and $L_3(\alpha,\beta)=-1$, (4.10) is
impossible. Thus, the lemma is proved.

Using the same method as in the proof of Lemma 4.3, we can obtain
the following lemma.

{\bf Lemma 4.4.} The equation $$15x^2+49^{2y}=2^{z+2},\ \ x,y,z\in
\mathbb{N}  $$ has no solution $(x,m,n)$.

{\bf Lemma 4.5.} If $D_1$ and $D_2$ satisfy (3.7), then
$N_2(D_1,D_2)=0$.

{\bf Proof.} Under the assumption, we suppose that (1.2) has a
solution $(x,m,n)$ with $2|m$. Since $2\nmid D_1D_2x$ and $D_2^m
\equiv 1(\rm mod \ 8)$, we have $D_1\equiv D_1x^2\equiv
2^{n+2}-D_2^m\equiv 7(\rm mod \ 8)$. Hence, by (3.7), we get
$\lambda=1$ and $$D_1X_1^2=2^{Z_1}-1,\  D_2=3.2^{Z_1}+1,\  Z_1\geq
3. \eqno{(4.11)}$$ By Lemmas 4.3 and 4.4, the lemma is true for
$Z_1 \in \{3,4 \}$, We may therefore assume that $Z_1\geq 5$.

Since $2|m$, we see from (1.2) that equation
$${X^\prime}^2+D_1{Y^\prime}^2=2^{Z^\prime+2},\ X^\prime,Y^\prime,Z^\prime \in \mathbb{Z},\ \textrm{gcd}(X^\prime,Y^\prime)=1,\ Z^\prime>0 \eqno{(4.12)}$$
has the solution
$$(X^\prime,Y^\prime,Z^\prime)=(D_2^{m/2},x,n).\eqno{(4.13)}$$ Apply
Lemma 2.13 to (4.13), we have $$n=Z_1^{\prime}t,t\in \mathbb{N},
\eqno{(4.14)}$$
$$\frac{D_2^{m/2}+x\sqrt{-D_1}}{2}=\lambda_1{\left(\frac{X_1^\prime+\lambda_2Y_1^\prime \sqrt{-D_1}}{2}  \right)}^t, \lambda_1,\lambda_2 \in \{\pm 1\}, \eqno{(4.15)}$$
where $(X_1^\prime,Y_1^\prime,Z_1^\prime)$ is the least solution of
(4.12).

Since $1^2+D_1X_1^2=2^{Z_1}$, by Lemma 2.17, the least solution
$(X_1^\prime,Y_1^\prime,Z_1^\prime)$ of (4.12) satisfies either
$$(X_1^\prime,Y_1^\prime,Z_1^\prime)=(1,X_1,Z_1-2) \eqno{(4.16)}$$
or $${X_1^\prime}^2-D_1{Y_1^\prime}^2=2\lambda, X_1^\prime
Y_1^\prime=X_1, Z_1^\prime=\frac{1}{2}(Z_1-2), \lambda \in \{\pm
1\}. \eqno{(4.17)}$$

We first consider the case (4.16), then (4.14) and (4.15) can be
written as $$n=(Z_1-2)t, t \in \mathbb{N}, \eqno{(4.18)}$$
$$\frac{D_2^{m/2}+x\sqrt{-D_1}}{2}=\lambda_1{\left(\frac{1+\lambda_2X_1 \sqrt{-D_1}}{2}  \right)}^t, \lambda_1,\lambda_2 \in \{\pm 1\}. \eqno{(4.19)}$$
Form (4.19), we get $X_1|x$, and hence, we have $x=X_1y$, where $y
\in \mathbb{N}$. Substituting it into (1.2), by (4.11) and (4.18),
we get $$D_1x^2+D_2^m=(2^{Z_1}-1)y^2+D_2^m=2^{n+2}=2^{(Z_1-2)t+2}.
\eqno{(4.20)}$$

Let $$\alpha=\frac{1+X_1\sqrt{-D_1}}{2},\ \
\beta=\frac{1-X_1\sqrt{-D_1}}{2}. \eqno{(4.21)}$$ By (4.11), (4.19)
and (4.21), we have
$$D_2^{m/2}={(3.2^{Z_1}+1)}^{m/2}=\lambda_1(\alpha^t+\beta^t).\eqno{(4.22)}$$
Since $\alpha+\beta=1$ and $\alpha\beta=2^{Z_1-2}$, by Lemma 2.1,
we have $$\alpha^t+\beta^t=\sum_{i=0}^{[t/
2]}(-1)^i\left[^t_i\right]2^{(Z_1-2)i}.\eqno{(4.23)}$$ Since
$\binom{t}{0}=1$, $\binom{t}{1}=t$ and $Z_1\geq 5$, compare (4.22)
and (4.23) we obtain $\lambda_1=1$ and
\begin{displaymath}
t\equiv \left\{\begin{array}{ll}
4\ (\rm mod\ 8), \ \textrm{if 2$\parallel m$};\\
0\ (\rm mod\ 8), \ \textrm{if 4$| m$}.
\end{array} \right.\eqno{(4.24)}
\end{displaymath}

Since $4|t$ by (4.24), we have $t=4s$, where $s\in \mathbb{N}$.
Hence, by (4.20), we get
$$2^{2(Z_1-2)s+1}+D_2^{m/2}=Af^2, 2^{2(Z_1-2)s+1}-D_2^{m/2}=Bg^2, \eqno{(4.25)}$$
$$AB=2^{Z_1}-1=D_1X_1^2,A,B,f,g\in \mathbb{N}, \textrm{gcd}(A,B)=\textrm{gcd}(f,g)=1.$$
From (4.25), we have $$Af^2+Bg^2=2^{2(Z_1-2)s+2}. \eqno{(4.26)}$$ By
Lemma 2.13, we see from (4.16) that the equation
$${X^{\prime\prime}}^{2}+(2^{Z_1}-1){Y^{\prime\prime}}^{2}=2^{Z^{\prime\prime}+2},X^{\prime\prime},
Y^{\prime\prime},Z^{\prime\prime}\in
\mathbb{Z},\textrm{gcd}(X^{\prime\prime},Y^{\prime\prime})=1,Z^{\prime\prime}>0
\eqno{(4.27)}$$ is bound to have a solution $(X^{\prime\prime},
Y^{\prime\prime},Z^{\prime\prime})$ with $Z^{\prime\prime}=2(Z_1-2)s.$\\
Therefore, since $AB=2^{Z_1}-1$, by Lemma 2.16, we get from (4.26)
that
$$(A,B)=(1,2^{Z_1}-1) \ \ \textrm{or}\ \  (2^{Z_1}-1,1).\eqno{(4.28)}$$

On the other hand, by (4.25), we have
$$Af^2-Bg^2=2D_2^{m/2}=2(3\cdot2^{Z_1}+1)^{m/2}. \eqno{(4.29)}$$ Since
$2^{Z_1}\equiv ({\rm mod} \ AB)$ and $D_2\equiv
3\cdot2^{Z_1}+1\equiv 4(\rm mod\ 2^{Z_1}-1)$, we have
$(D_2/A)=(D_2/B)=1$, where $(\ast/\ast)$ is the Jacobi symbol.
Hence, by (4.25), we get
$$1=\left(\frac{-2B}{A}\right)=\left(\frac{-2}{A}\right), \
1=\left(\frac{2A}{B}\right)=\left(\frac{2}{B}\right).
\eqno{(4.30)}$$ The combination of (4.28) and (4.30) yields
$(A,B)=(1,2^{Z_1}-1)$. Substituting it into the first equality of
(4.25), we have $$f^2-D_2^{m/2}=2^{2({Z_1-2})s+1}. \eqno{(4.31)}$$

If $4|m$, then from (4.31) we get $f+D_2^{m/4}=2^{2(Z_1-2)s}$ and
$f-D_2^{m/4}=2$, whence we obtain
$$D_2^{m/4}=2^{2(Z_1-2)s-1}-1.\eqno{(4.32)}$$
Since $2(Z_1-2)s-1\geq 5$, apply Lemma 2.8 to (4.32), we get $m=4$
and $D_2= 3\cdot2^{Z_1}+1= 2^{2(Z_1-2)s-1}-1$. But, since $Z_1\geq
5$, it is impossible. So we have $2\parallel m$.

Since $2\parallel m$, we see from (4.24) that $4\parallel t$ and
$t=4s$, where $s\in \mathbb{N}$ with $2\nmid s$. Since
$\alpha^4+\beta^4=1-2^{Z_1}+2^{2Z_1-3}$ by (4.21), apply Lemma 2.1
to (4.22), we have $\lambda_1=1$ and
$$D_2^{m/2}=(3\cdot2^{Z_1}+1)^{m/2}=\alpha^t+\beta^t=
(\alpha^4+\beta^4)\left(\frac{(\alpha^4)^s+(\beta^4)^s}{\alpha^4+\beta^4}
\right) $$
$$=(1-2^{Z_1}+2^{2Z_1-3})\sum_{j=0}^{(s-1)/2}(-1)^{j}\left[^s_j\right]{(1-2^{Z_1}+2^{2Z_1-3})}^{s-2j-1}2^{4(Z_1-2)j}.\eqno{(4.33)}$$
By (4.33), we get
$2^{2Z_1-3}-2^{Z_1}+1|{(3\cdot2^{Z_1}+1)}^{m/2}$. Let
$d=\textrm{gcd}(2^{2Z_1-3}-2^{Z_1}+1,3\cdot2^{Z_1}+1)$. Then we
have $d|97$. Hence, we get
$$2^{2Z_1-3}-2^{Z_1}+1=97^k,k\in \mathbb{N}. \eqno{(4.34)}$$
Apply Lemma 4.2 to (4.34), we obtain $Z_1=5$ and $k=1$. It implies that $(D_1,D_2)=(31,97)$. Since (1.2) has a solution $(x,m,n)=(15,2,12)$ for $(D_1,D_2)=(31,97)$,
 by Lemma 3.4, we get $N_2(31,97)=1$.
Moreover, if $Z_1>5$, then (1.2) has no solution $(x,m,n)$ with
$2|m$ for the case (4.16).

We next consider the case (4.17). If $2\nmid t$ for (4.14), then
form (4.15) we get
$X_1^\prime|D_2^{m/2}={(3\cdot2^{Z_1}+1)}^{m/2}$. Since
$X_1^\prime|X_1$, $X_1|2^{Z_1}-1$ and
$\textrm{gcd}(2^{Z_1}-1,3\cdot2^{Z_1}+1)=1$, we get
$X_1^\prime=1$. Substituting it into (4.17), we obtain
$Y_1^\prime=1$ and $D_1=3\neq 2^{Z_1}-1$, a contradiction. So we
have $2|t$. Let $t=2t^\prime$, where $t^\prime \in \mathbb{N}$.
Since
$${\left(\frac{X_1^\prime+\lambda_2Y_1^\prime \sqrt{-D_1}}{2}
\right)}^{2}=\lambda \frac{1+\lambda\lambda_2X_1\sqrt{-D_1}}{2}$$ by
(4.17), we get from (4.14) and (4.15) that $$n=(Z_1-2)t^\prime,
t^\prime \in \mathbb{N}, \eqno{(4.35)}$$
$$\frac{D_2^{m/2}+x\sqrt{-D_1}}{2}=\lambda_1\lambda^{t^\prime}{\left(\frac{1+\lambda\lambda_2X_1\sqrt{-D_1}}{2} \right)}^{t^\prime},\lambda,\lambda_1,\lambda_2 \in \{\pm 1 \}. \eqno{(4.36)} $$

Obviously, (4.35) and (4.36) are equal to (4.18) and (4.19)
respectively. Thus, by the conclusion of case (4.16), the lemma is
proved.

{\bf Proof of Theorem B.} Let $(x,m,n)$ be a solution of (1.2). Then
we have $D_1\equiv 7(\rm mod\ 8)$ if $2\nmid m$. Therefore, by
Lemmas 3.2 and 4.5, if $N_1(D_1,D_2)>0$, then $N_1(D_1,D_2)\leqslant
2$ except for $N(3,5)=N(5,3)=4$ and $N(13,3)=N(31,97)=3$. Further,
by Lemmas 3.1, 3.4 and 4.5, if $N_2(D_1,D_2)>0$, then
$N(D_1,D_2)\leqslant 2$. Thus, the theorem is proved.

\section*{References}

[1] M. Abouzaid, Les nombres de Lucas et Lehmer sans diviseur
primitif, J. Th\'{e}orie Nombres Bordeaux, 2006, 18(2): 299-313.

[2]  Y. Bilu, G. Hanrot and P. M. Voutier (with an appendix by M.
Mignotte), Existence of primitive divisors of Lucas and Lehmer
numbers, J. Reine Angew. Math., 2001, 539: 75-122.

[3] G. D. Birkhoff and H. S. Vandiver, On the integral divisors of
$a^n-b^n$, Ann. of Math. (2), 1904, 5: 173-180.

[4] Y. Bugeaud, On some exponential diophantine equations, Monatst.
Math., 2001, 132(1): 93-97.

[5] Y. Bugeaud and M. Mignotte, On the diophantine equation
$(x^n-1)/(x-1)=y^q$ with negative $x$, In: M. A. Bennett, Number
theory for the millennium I, Urbana-Champaign, IL, 2002, 145-151.

[6] Y. Bugeaud, M. Mignotte and S. Siksek, Classical and modular
approaches to exponential diophantine equation I: Fibonacci and
Lucas perfect powers, Ann. of Math. (2), 2006, 163(2): 969-1018.

[7] Y. Bugeaud and T. N. Shorey, On the number of solutions of the
generalized Ramanujan-Nagell equation, J. Reine Angew. Math., 2001,
539: 55-74.

[8] M.-H. Le, On the diophantine equation  $D_1x^2+D_2=2^{n+2}$,
Acta Arith., 1993, 64(1): 29-41.

[9] M.-H. Le, On the diophantine equation $x^2+D^m=2^{n+2}$,
Commen. Math. Univ. St. Pauli, 1994, 43(2): 127-133.

[10] R. Lidl and H. Niederreiter, Finite fields, Addison - Wesley,
Reading, MA, 1983.

[11] W. Ljunggren, Oppgave nr 2, Norsk. Mat. Tidsskrift, 1943,
25(1): 29.

[12]  W. Ljunggren, Noen setninger om ubestemte likninger av
formen $(x^n-1)/(x-1)=y^q$, Norsk. Mat. Tidsskrift, 1943, 25(1):
17-20.

[13] P. Mih$\check{a}$ilescu, Primary cyclotomic units and a proof
of Catalan's conjecture, J. Reine Angew. Math., 2004, 572: 167-195.

[14] L. J. Mordell, Diophantine equations, Academic Press, London,
1969.

[15] T. Nagell, Des \'{e}quations ind\'{e}termin\'{e}es
$x^2+x+1=y^n$ et $x^2+x+1=3y^n$, Norsk. Mat. Forenings Skr. Series
I, 1921, 2(1): 12-14.

[16] T. Nagell, \"{U}ber die rationaler punkte auf einigen
kubischen kurven, Tohoku Math. J., 1924, 24(1): 48-53.

[17] T. Nagell, L$\phi$sning til oppgave nr 2, Norsk. Mat.
Tidsskrift, 1948, 30(1): 62-64.

[18] S. Ramanujan, Question 464, J. Indian Math. Soc., 1913, 5(2):
120.

[19] P. M. Voutier, Primitive divisors of Lucas and Lehmer
sequences, Math. Comp., 1995, 64(5): 869-888.

\newpage

\end{CJK*}
\end{document}